\documentclass[reqno,12pt]{amsart}

\usepackage[latin1]{inputenc}
\usepackage{amsmath,amsthm,amssymb,amsfonts}
\usepackage[all]{xy}

\newcommand{\R}{\mathbb{R}}
\newcommand{\Z}{\mathbb{Z}}
\newcommand{\N}{\mathbb{N}}
\newcommand{\U}{\mathcal{U}}
\newcommand{\V}{\mathcal{V}}
\newcommand{\W}{\mathcal{W}}
\newcommand{\del}{\mathcal{D}}

\newcommand{\back}{\!\!\!\!\!\!\!\!}
\newcommand{\id}{\mathrm{id}}
\newcommand{\lmap}{\ell}
\newcommand{\geo}{\boldsymbol|}
\newcommand{\st}{\mathrm{st}}
\newcommand{\aw}{\mathrm{AW}}

\def\paritem#1{%
  \smallskip
  \setbox0=\hbox{#1\enspace}
  \par\noindent
  \ifnum\wd0>\parindent\box0
  \else\hbox to\parindent{\box0\hfill}\fi\ignorespaces}

\theoremstyle{plain}
\newtheorem{theorem}{Theorem}[section]
\newtheorem{proposition}[theorem]{Proposition}
\newtheorem{lemma}[theorem]{Lemma}
\newtheorem{corollary}[theorem]{Corollary} 

\theoremstyle{definition}
\newtheorem{defi}[theorem]{Definition}
\newtheorem{eks}[theorem]{Example}
\newtheorem{rem}[theorem]{Remark}
\newtheorem{remarks}[theorem]{Remarks}

\numberwithin{equation}{section}
\makeatletter
\let\c@equation\c@theorem
\makeatother

\title{Integration of simplicial forms and Deligne cohomology}
\author[J.~L~Dupont]{Johan L.~Dupont}
\address{Department of Mathematics   \\
University of Aarhus \\
DK-8000 {\AA}rhus, Denmark}
\email[J.~L.~Dupont]{dupont@imf.au.dk}
\author[R.~Ljungmann]{Rune Ljungmann}
\address{Department of Mathematics   \\
University of Aarhus \\
DK-8000 {\AA}rhus, Denmark}
\email[R.~Ljungmann]{runel@imf.au.dk}

\begin{document}
\thanks{
Work supported in part by the Danish Natural Science Research Council and the European Union Network EDGE} 
\date{\today}
\begin{abstract}
We present two approaches to constructing an integration map for smooth
Deligne cohomology. The first is defined in the simplicial model, where a class in
Deligne cohomology is represented by a simplicial form, and the second in a
related but more combinatorial model.
\end{abstract}

\maketitle
\tableofcontents
\section{Introduction}
For the construction of invariants for families of bundles, integration
along the fiber is usually applied in order to obtain forms defined on the
parameter space. In the case of families of bundles with connection the
classical Chern-Weil theory gives rise to invariants living in smooth
Deligne cohomology, and hence a notion of integration along the fiber is  needed in this setting (see e.g. Freed \cite{F} or Dupont-Kamber \cite{DK}). In this paper, we introduce two different constructions of this map. The first one is defined in the simplicial model for smooth Deligne cohomology introduced in \cite{DK}, where a class in the smooth Deligne cohomology $H^{l+1}_\del(Z,\Z)$ is represented by a simplicial form $\omega\in\Omega^l(|N\U|)$, for $\U$ a covering of $Z$. We prove that:
\begin{theorem}\label{firstmain}
Given a fiber bundle $\pi:Y\to Z$ with compact, oriented $n$-dimensional fibers and suitable coverings $\V$ and $\U$. Then there is a map
\[\int_{[Y/Z]}:\Omega^{*+n}(|N\V|)\to \Omega^*(|N\U|).\]
It satisfies a Stokes' formula, and thus if $\partial Y=\emptyset$ induces a map
\[\pi_!:H^{*+n}_\del(Y,\Z)\to H^*_\del(Z,\Z)\]
in smooth Deligne cohomology independent of all choices.
\end{theorem}
The second construction defines the map in a more combinatorial model where the cohomology classes are represented by simplicial forms living in the 'triangulated nerve' $|NK|$ associated to a triangulation $|K|\to|L|$ of the bundle. This allows us to state the following useful theorem in the case where the fiber has boundary:
\begin{theorem}\label{secondmain}
Assume that $\partial Y\neq \emptyset$ then for a form $\omega\in\Omega^{*+n}(|N\V|)$ representing an element in smooth Deligne cohomology, the form
\[\int_{K/L}\omega\in \Omega^*(|NL|)/d\Omega^{*-1}(|NL|)\]
depends only on the triangulation of $\partial Y\to Z$.
\end{theorem}

There are other approaches to the subject in the literature. In Hopkins-Singer \cite{HS}, a cochain model for the Cheeger-Simons differential characters is given and an integration map is constructed by embedding the bundle in a larger trivial one. In the \v Cech-de Rham model, Gomi-Terashima \cite{GT} have introduced a combinatorial formula that uses a triangulation of the fiber. Unfortunately their formula is given for product bundles only, and it is not immediately clear how to generalise it to the case of a general fiber bundle. We hope to demonstrate that the approach with simplicial forms is a natural generalisation of the usual integration map.

We will start by giving a short description of smooth Deligne cohomology
both in the usual \v Cech-de Rham model and in the simplicial model
introduced in \cite{DK} in \S2. In \S3, we introduce the concept of
\emph{prism complexes} which is a generalisation of simplicial sets well suited for fiber bundles. It will provide a convenient framework for the
constructions in \S\S4-5. In \S4, we construct an integration map in the
simplicial model by choosing suitable coverings of the fiber bundle and a
set of partitions of unity. In \S5, we introduce a more combinatorial model
closely related to the simplicial approach. By using an Alexander-Whitney
type map, we then give a combinatorial integration formula. Finally the two
approaches are shown to induce the same map in smooth Deligne
cohomology.

\textbf{Acknowledgment} The authors would like to thank Franz Kamber, Ulrich Bunke and Marcel B\" okstedt for useful discussions during the preparation of this paper.


\section{Smooth Deligne cohomology}
Here follows a short introduction to smooth Deligne cohomology. The Deligne cohomology groups are usually constructed as the hypercohomology of a certain sequence of sheaves. We will however restrict ourselves to the corresponding concrete \v Cech description. For a more comprehensive exposition see Brylinski \cite{B}.

Let $Z$ be a smooth manifold of dimension $m$ and let $\U=\{U_i\}_{i\in I}$ be a 'good' open cover of $Z$. (That is every non-empty intersection of sets from the covering is contractible).

Let $\check{\Omega}^{p,q}(\U)=\check{C}^p(\U,\underline{\Omega}^q)$ be the ordinary \v Cech-de Rham complex. It is well-known that the chain-map
\[\varepsilon^*:\Omega^q(Z)\to \check{\Omega}^{0,q}(\U),\]
induced by the inclusion $\varepsilon:\sqcup U_i\to Z$ gives an isomorphism
\[H_{\mathrm{dR}}^*(Z)\to H^*(\check{\Omega}^*(\U))\]
in homology. We also have an inclusion of the ordinary \v Cech-complex with integer coefficients
\[\check{C}^p(\U,\Z)\to \check{\Omega}^{p,0}(\U)\]
which gives us the quotient complex
\[\check{\Omega}_{\R/\Z}^*(\U)=\check{\Omega}^*(\U)/\check{C}^*(\U,\Z).\]
\begin{defi}
\paritem{$1.$}
An (Hermitian line) \emph{$l-1$-gerbe} on $Z$ is an $l$-cocycle in $\check{C}^l(\U,\underline{\R/\Z})$ or equivalently a $\theta\in\check\Omega_{\R/\Z}^{l,0}(\U)$ with $\delta\theta=0$.
\paritem{$2.$}
A \emph{connection} $\omega$ in an $l-1$-gerbe $\theta$ is a $\omega=(\omega_0,\dots,\omega_l)\in\check\Omega_\R^l(\U)$, where $\omega_i\in\check\Omega_\R^{i,l-i}(\U)$, so that $\omega_l\equiv -\theta \mod \Z$ and $\omega$ is a cycle in $\check\Omega_{\R/\Z}^*/\varepsilon^*\Omega^*(Z)$. 
\paritem{$3.$}
Two $l-1$-gerbes $\theta$ and $\theta'$ with connections $\omega$ and $\omega'$ are equivalent if $\omega$ and $\omega'$ are cohomologous in $\check\Omega_{\R/\Z}^*(\U)$. The set of equivalence classes $[\theta,\omega]$ is denoted $H^{l+1}_\del(Z,\Z)$ and is called the $l+1$'st \emph{smooth Deligne cohomology} group.
\end{defi}
\begin{remarks}
\paritem{$1.$}
Note that $H_\del^{l+1}(Z,\Z)$ is the cohomology of the sequence
\[\check\Omega_{\R/\Z}^{l-1}(\U)\stackrel{d}{\to} \check\Omega_{\R/\Z}^l(\U)\stackrel{d}{\to}\check\Omega_{\R/\Z}^{l+1}(\U)/\varepsilon^*\Omega^{l+1}(Z).\]
\paritem{$2.$}
That $\omega$ is a cycle in $\check\Omega_{\R/\Z}^*/\varepsilon^*\omega^*(Z)$ is equivalent to the relations 
\[\delta\omega_{i-1}+(-1)^id\omega_i=0,\quad i=0,\dots,l\]
and 
\[\delta\omega_l\equiv 0 \mod \Z.\]
\paritem{$3.$}
Our definition of a gerbe is to some extend an abuse of language, since a gerbe is actually a well-defined geometrical object, so that the set of isomorphism classes of gerbes (with band $\R/\Z$) over $Z$ is isomorphic to $H^2(Z,\underline{\R/\Z})$, this corresponds to our case $l=2$. Our viewpoint is equivalent to identifying a line bundle with its defining cocycle. If the reader finds this inconvenient, he/she can simply choose to substitute 'gerbe' with 'gerbe data'.
\end{remarks}
\begin{proposition}
\paritem{$1.$}
We have a commutative diagram
\begin{eqnarray*}
\xymatrix{
H_\del^{l+1}(Z,\Z) \ar[d]_{\delta_*} \ar[r]^{d_*} & \Omega_{\mathrm{cl}}^{l+1}(Z) \ar[d]_{I}\\
H^{l+1}(Z,\Z) \ar[r]& H^{l+1}(Z,\R)
}
\end{eqnarray*}
where $\Omega_{\mathrm{cl}}^{l+1}(Z)$ is the set of closed $l+1$-forms with integral periods.
\paritem{$2.$}
There is a short exact sequence
\begin{equation}\label{shortexactdeligne}
0\to H^l(Z,\R/\Z)\to H^{l+1}_\del(Z,\Z)\stackrel{d_*}{\to} \Omega_{\mathrm{cl}}^{l+1}(Z)\to 0,
\end{equation}
\end{proposition}
\begin{proof}
\paritem{$1.$}
Note that since $\delta d\omega_0=d\delta\omega_0=d^2\omega_1=0$ then $F_\omega=d\omega_0$ is actually a globally defined, closed $l+1$-form. $F_\omega$ is called the \emph{curvature} of $\omega$. $d_*$ is the map sending $\omega$ to $F_\omega$. $\delta_*$ is just the connecting homomorphism for the short exact sequence
\[0\to \Z\to \underline \R \to \underline{\R/\Z}\to 0\]
and $I$ is the de Rham map. Now, commutativity of the diagram follows from the fact that $d\omega_0-\delta\theta=d\omega_0+\delta\omega_l=D\omega$ in $\check\Omega_\R^*(\U)$.
\paritem{$2.$}
The kernel of $d_*$  is the 'gerbes with \emph{flat} connection'. Since a gerbe with flat connection is actually a cocycle in the full complex $\check{\Omega}_{\R/\Z}^*(\U)$, we see that the kernel is in fact $H^l(Z,\R/\Z)$.
\end{proof}

There is also a description of gerbes with connection in terms of the differential characters of Cheeger-Simons \cite{CS}. Indeed there is an explicit isomorphism $H^l_\del(Z,\Z)\cong \hat H^l(Z,\Z)$ given in e.g. Dupont-Kamber \cite{DK}.


\subsection{Simplicial forms}\label{afsnit:simplicielleformer}
In \cite{DK}, there is  given a description of gerbes with connection in terms of simplicial forms, which we will briefly recall. Given an open cover $\U=\{U_i\}$ of $Z$ we have the nerve $N\U=\{N\U(p)\}$ of the covering where 
\[N\U(p)=\bigsqcup_{i_0,\dots,i_p} U_{i_0}\cap\cdots\cap U_{i_p}.\]
We denote $U_{i_0}\cap\cdots\cap U_{i_p}$ by $U_{i_0\dots i_p}$ in the following.

$N\U$ is a simplicial manifold where the face maps
\[d_j:U_{i_0\dots i_p}\to U_{i_0\dots \hat{i_j}\dots i_p}\]
and degeneracy maps
\[s_j:U_{i_0\dots i_p}\to U_{i_0\dots i_j i_j\dots i_p}\]
 are just inclusions.
\begin{defi}
A \emph{simplicial $n$-form} $\omega=\{\omega^{(p)}\}$ on $N\U$ consists of forms $\omega^{(p)}\in\Omega^n(\Delta^p\times N\U(p))$ which satisfy the relations
\[(\varepsilon_j\times \id)^*\omega^{(p)}=(\id\times d_j)^*\omega^{(p-1)},\]
where $\varepsilon_j:\Delta^{p-1}\to\Delta^p$ denotes the ordinary $j$'th face map. We denote the set of simplicial forms on $N\U$ by $\Omega^*(||N\U||)$. If the forms also satisfy the relations
\[(\eta_j\times \id)^*\omega^{(p-1)}=(\id\times s_j)^*\omega^{(p)},\]
where $\eta_j:\Delta^p\to\Delta^{p-1}$ is the ordinary $j$'th degeneracy map, the forms are called \emph{normal}. The set of normal forms is denoted $\Omega^*(|N\U|)$.
\end{defi}
\begin{rem}\label{rempermutation}
Our index sets will always be assumed to be ordered and it is then customary to consider only ordered $(p+1)$-tuples, that is for a tuple $(i_0,\dots,i_p)$ we have $i_0\leq \dots \leq i_p$. Later when we move on to prism complexes this will in some instances be annoying. Instead we demand that for a permutation $\sigma \in\Sigma(p)$ the normal forms also satisfy the relation
\[\tilde\sigma^*\omega=\omega\]
where $\tilde\sigma:\Delta^p\times U_{i_0\dots i_p}\to \Delta^p\times U_{i_{\sigma(0)}\dots i_{\sigma(p)}}$ on the first factor is the simplicial map that permutes the vertices of $\Delta^p$ according to $\sigma$  and on the second factor is the identity.
\end{rem}

We have a direct sum decomposition
\[\Omega^n(|N\U|)=\bigoplus_{p+q=n} \Omega^{p,q}(|N\U|)\]
where $\Omega^{p,q}(|N\U|)$ is the set of forms that are of degree $p$ in the barycentric coordinates on the simplex in the product $\Delta^k\times N\U(k)$ for $k\geq p$.

There is a chain map
\[I_\Delta:\Omega^{p,q}(|N\U|)\to \check\Omega^{p,q}(\U)\]
given by $I_\Delta(\omega)=\int_{\Delta^p}\omega^{(p)}$.
This map gives an isomorphism in homology. In fact it has a right inverse given on $\Delta^k\times N\U(k)$ by
\[E(\omega)=p!\sum_{|I|=p}\omega_I\wedge d_I^*\omega\]
Where $I=(i_0,\dots,i_p)$ is a sequence of integers $0\leq i_0\leq \cdots\leq i_p\leq k$, $\omega_I=\sum_{j=0}^p (-1)^jt_{i_j}dt_{i_0}\wedge \hat{dt_{i_j}}\wedge dt_{i_p}$ are the \emph{elementary forms} on $\Delta^k$ and $d_I:N\U(k)\to N\U(p)$ is $d_I=d_{j_1}\cdots d_{j_l}$ where $0\leq j_l\leq \cdots\leq j_1\leq k$ is the complementary sequence of $I$ (see Dupont \cite{D} for details).

The natural inclusion $\sqcup U_i \to Z$ also induces a map
\[\varepsilon'^*:\Omega^*(Z) \to \Omega^*(|N\U|),\]
so we get the following commutative diagram of homology isomorphisms:
\begin{eqnarray*}
\xymatrix{
\Omega^n(Z) \ar[dr]_{\varepsilon^*} \ar[r]^{\varepsilon'^*} & \Omega^n(|N\U|) \ar[d]_{I_{\Delta}}\\
& \check\Omega^n(\U)
}
\end{eqnarray*}

We need a notion of integral simplicial forms in order to imitate the construction in the previous section. 
\begin{defi}
A form $\omega\in\Omega^*(|N\U|)$ is called \emph{discrete} if it is constant
with respect to any local coordinates on the nerve. Furthermore it is called
\emph{integral} if $I_\Delta(\omega)\in \check C^*(\U,\Z)$. The chain complex of integral forms is denoted $\Omega^*_{\Z}(|N\U|)$
\end{defi}

\begin{proposition}
We have the following isomorphisms
\paritem{$1.$}
\[H^n(\Omega^*_\Z(|N\U|))\cong H^n(C^*(\U,\Z))=H^n(Z,\Z).\]
\paritem{$2.$}
If we define 
\[\Omega^*_{\R/\Z}(|N\U|)=\Omega^*(|N\U|)/\Omega^*_\Z(|N\U|)\]
then also
\[H^n(\Omega^*_{\R/\Z}(|N\U|))\cong H^n(\check\Omega^*_{\R/\Z}(\U))\cong H^n(Z,\R/\Z).\]
\paritem{$3.$}
$I_\Delta$ induces an isomorphism from the cohomology of the sequence
\begin{equation}\label{simplicialcohom}
\Omega^{l-1}_{\R/\Z}(|N\U|)\stackrel{d}{\to}\Omega^l_{\R/\Z}(|N\U|)\stackrel{d}{\to}\Omega^{l+1}_{\R/\Z}(|N\U|)/\varepsilon^*\Omega^{l+1}(Z)
\end{equation}
to $H^{l+1}_\del(Z,\Z)$.
\end{proposition}
\begin{proof}
\paritem{$1.$}
The map $I_\Delta$ takes integral forms to integral cochains by definition. It induces an isomorphism in cohomology since the map $E$ takes integral cochains to integral forms and the chain homotopies from $\id$ to $E\circ I_\Delta$ given in \cite{D} are easily seen to map integral forms to integral forms. 
\paritem{$2.$} 
Follows easily from the above.
\paritem{$3.$} 
Since the cohomology group of (\ref{simplicialcohom}) fits into the same short exact sequence as in (\ref{shortexactdeligne}) the $5$-lemma gives us that $I_\Delta$ is an isomorphism.
\end{proof}
\begin{corollary}
Every class in $H_\del^{l+1}(Z,\Z)$ can be represented by an $l-1$-gerbe $\theta$ with connection $\omega$, where $\omega=I_\Delta(\Lambda)$ for some $\Lambda\in\Omega^l(|N\U|)$ and
\[d\Lambda=\varepsilon^*\alpha-\beta,\quad \alpha\in\Omega^{l+1}(Z),\quad \beta\in\Omega_\Z^{l+1}(|N\U|).\]
\end{corollary}

\section{Prism complexes}
The notion of a 'prism complex' and 'prismatic' decomposition has occurred (implicitly or explicitly) in many different contexts. We refer to Akyar-Dupont-Ljungmann \cite{ADL} for further details and references.

A prism complex is a generalisation of a simplicial set (or manifold) well suited for fiber bundles.
A \emph{prism complex} $P=\{P_p\}$ is a collection of $p+1$-simplicial sets
$P_p$. That is for each set of positive integers $(q_0,\dots,q_p)$ we have
sets $P_{p,q_0\dots q_p}$ with face and degeneracy maps
$d^i_j:P_{p,q_0\dots q_p} \to P_{p,q_0\dots q_i-1\dots q_p}$ and
$s^i_j:P_{p,q_0,\dots,q_p}\to P_{p,q_0\dots q_i+1\dots q_p}$,
$i=0,\dots,p$, $j=0,\dots q_i$ satisfying the usual relations. Furthermore
we want another set of simplicial (i.e. commuting with the $d^i_j$'s
and $s^i_j$'s) face maps $d_i:P_{p,q_0\dots q_p}\to P_{p-1,q_0\dots \hat{q_i}\dots q_p}$ and degeneracy maps $s_i:P_{p,q_0\dots q_p}\to P_{p+1,q_0,\dots q_iq_i\dots q_p}$ so that $(P_p,d_i,s_i)$ becomes an ordinary simplicial set. Note that in some applications the last set of degeneracy maps does not exist naturally so in these cases $(P_p,d_i)$ is only a $\Delta$-set. 
As with ordinary simplicial sets we can for each $p$ form the geometric and fat realisations $|P_p|$ and $||P_p||$, that is, the quotients of 
\[\bigsqcup_{q_0\dots q_p} \Delta^{q_0}\times \cdots \times \Delta^{q_p}\times P_{p,q_0\dots q_p}\]
where we divide out by the equivalence relations generated by the face and degeneracy maps 
\[\varepsilon^i_j:\Delta^{q_0\dots q_i\dots q_p}\to\Delta^{q_0\dots q_i+1\dots q_p}\]
and (in case of the geometric realisation)
\[\eta^i_j:\Delta^{q_0\dots q_i\dots q_p}\to \Delta^{q_0\dots q_i-1\dots q_p}\]
(where $\Delta^{q_0\dots q_p}$ is short hand notation for the prism $\Delta^{q_0}\times\cdots\times \Delta^{q_p}$).
 
The face and degeneracy maps $d_i$ and $s_i$ now induce a  structure of a simplicial set on $|P_p|$ ($||P_p||$) by acting as the projection and the diagonal on $\Delta^{q_0}\times\cdots\times \Delta^{q_p}$ respectively. That is let $\pi_i:\Delta^{q_0\dots q_p}\to\Delta^{q_0\dots \hat{q_i}\dots q_p}$ be the projection that deletes the $i$'th coordinate and let $\Delta_i:\Delta^{q_0\dots q_p}\to \Delta^{q_0\dots q_iq_i\dots q_p}$ be the diagonal map that repeats the $i$'th factor. Then we can form the geometric realisation
\[\geo P.\geo=\bigsqcup_{p\geq 0}\Delta^p\times |P_p|/\sim\]
where the equivalence relation is generated by
\[(\varepsilon_i t,s,x)\sim (t,\pi_is,d_ix),\quad t\in\Delta^{p-1},\,s\in \Delta^{q_0\dots q_p},\,x\in P_{p,q_0\dots q_p}\]
and
\[(\eta_i t,s,x)\sim (t,\Delta_is,s_ix),\quad t\in\Delta^{p+1},\,s\in \Delta^{q_0\dots q_p},\,x\in P_{p,q_0\dots q_p}\]

\begin{eks}\label{eksbundletriangulations}
Given a smooth fiber bundle $\pi:Y\to Z$ with $\dim Y=m+n$, $\dim Z=m$ and compact fibers, possibly with boundary, a theorem of Johnson \cite{J} gives us smooth triangulations $K$ and $L$ of $Y$ and $Z$ respectively and a simplicial map $\pi':K\to L$ so that the following diagram commutes
\begin{eqnarray*}
\xymatrix{
|K| \ar[d]^{|\pi'|} \ar[r]^{\cong} & Y \ar[d]_{\pi} \\
|L| \ar[r]^{\cong} & Z
}
\end{eqnarray*}
Here the horisontal maps are homeomorphisms which are smooth on each simplex. Furthermore given such a triangulation of $\partial Y \to Z$ we can also extend it to a triangulation of $Y\to Z$.

Now the geometric idea is that if $z\in Z$ lies in the interior of a $p$-simplex of $L$ then the fiber over $z$ is in a canonical way decomposed into $p+1$-fold prisms of the form $\Delta^{q_0\dots q_p}$ as above. Formally we define the prismatic complex $PS(K/L)$ by letting $PS_p(K/L)_{q_0\dots q_p}\subseteq S_{p+q_0+\cdots+q_p}(K)\times S_p(L)$ be the subset of simplices $(\tau,\eta)$ so that $q_i+1$ of the vertices in $\tau$ lies over the $i$'th vertex in $\eta$. Then we have face and degeneracy operators defined in the obvious way. In particular this gives us boundary maps in the 'fiber direction' of the associated chain complex
\[\partial^i_F:PC_p(K/L)_{q_0\dots q_p}\to PC_p(K/L)_{q_0\dots q_i-1\dots q_p}\]
defined by $\partial^i_F=\sum (-1)^jd^i_j$, ($\partial^i_F=0$ for $q_i=0$), and also a 'total' boundary map 'along the fiber'
\[\partial_F=\partial^0_F+(-1)^{q_0+1}\partial^1_F+\cdots+(-1)^{q_0+\cdots+q_{p-1}+p}\partial^p_F.\]
Also there is a 'horisontal' boundary map  
\[\partial_H=\partial_0+(-1)^{q_0+1}\partial_1+\cdots+(-1)^{q_0+\cdots+q_{p-1}+p}\partial_p,\]
where
\begin{equation*}
\partial_i=\left\{\begin{array}{ll}
0 &\textrm{if } q_i>0\\
d_i & \textrm{if }q_i=0
\end{array}\right.
\end{equation*}
so that $\partial=\partial_F+\partial_H$ is compatible with the chain map
in the ordinary chain complex $C^*(K)$ for $K$.

There is a natural 'prismatic triangulation' homeomorphism
\[\lmap: \geo PS(K/L)\geo\stackrel{\cong}{\to}|K|\]
induced by
\[\lmap(t,s^0,\dots,s^p,(\tau,\eta))=(t_0s^0,\dots,t_ps^p,\tau)\]
for $(t,s,\tau)\in \Delta^p\times \Delta^{q_0\dots q_p}\times PS_p(K/L)_{q_0\dots q_p}$.
Note that if $\stackrel{\circ}{\sigma}$ is an open $p$-simplex in $L$ then $\lmap$ provides a natural trivialisation of $|K|$ over $\stackrel{\circ}{\sigma}$
\[\stackrel{\circ}{\sigma}\times|PS_p(K/\sigma)|\stackrel{\cong}{\to}|K|_{|\sigma}\]

\end{eks}

\begin{eks}\label{eksprismaticnerve}
Another example in the category of manifolds, comes from the nerve of
compatible open coverings of the total space and the base space. That is,
given a covering $\U=\{U_i\}$ of $Z$ we have a covering
$\W=\{W_i=\pi^{-1}(U_i)\}$ of $Y$, and for each $i$, $\V^i$ is an open cover
of $W_i$. This gives a covering $\V=\cup \V^i$ of $Y$ (with lexicographically
ordered index set). Then we put
\[P_pN(\V/\U)_{q_0\dots q_p}=\bigsqcup V^{i_0}_{j^0_0}\cap\cdots \cap  V^{i_0}_{j^0_{q_0}}\cap\cdots\cap  V^{i_p}_{j^p_{q_p}}\]
with $V^i_j\in \V^i$, and face and degeneracy maps are inclusions similarly to the simplicial case in section \ref{afsnit:simplicielleformer}. In the following, we will denote  $V^{i_0}_{j^0_0}~\cap\cdots~\cap~V^{i_p}_{j^p_{q_p}}$ by $V_{j^0_0\dots j^p_{q_p}}$.

A useful special case of this situation occurs in the context of example
\ref{eksbundletriangulations} above with the coverings consisting of the (open) stars of the triangulations of $K$ and $L$. More precisely $\U=\{U_i=\st(a_i)\}$ where $a_i\in L^0$ is a $0$-simplex in $L$ and $\V^i=\{V^i_j=\st(b^i_j)\}$ where $b^i_j\in \pi^{-1}(a_i)\cap K^0$. Note that the discrete prismatic nerve of this covering is just $PS(K/L)$.
\end{eks}


\subsection{A de Rham theorem}
As a straightforward generalisation of simplicial forms, we introduce the complex of (normal) prismatic forms on the prism complex in the above example \ref{eksprismaticnerve}.
\begin{defi}
A \emph{prismatic $n$-form} is a collection $\omega=\{\omega_{q_0\dots q_p}\}$ of forms $\omega_{q_0\dots q_p}\in \Omega^n(\Delta^p\times\Delta^{q_0\dots q_p}\times P_pN\V/\U_{q_0\dots q_p})$ satisfying the relations 
\[(\id\times \varepsilon^i_j\times \id)^*\omega_{q_0\dots q_p}=(\id\times \id\times d^i_j)^*\omega_{q_0\dots q_i-1\dots q_p}\]
and
\[(\varepsilon_i\times \id\times \id)^*\omega_{q_0\dots q_p}=(\id\times \pi_i\times d_i)^*\omega_{q_0\dots \hat{q_i}\dots q_p}.\]
A form is called \emph{normal} if it also satisfies the relations
\[(\id\times\eta^i_j\times \id)^*\omega_{q_0\dots q_i-1\dots q_p}=(\id\times \id\times s^i_j)^*\omega_{q_0\dots q_p}\]
and
\[(\eta_i\times \id\times \id)^*\omega_{q_0\dots q_p}=(\id\times \Delta_i\times s_i)^*\omega_{q_0\dots q_i q_i\dots q_p}.\]
The complex of normal simplicial forms is denoted by
\[\Omega^*(\geo PN\V/\U\geo).\]
\end{defi}
As in the simplicial case we have a direct sum decomposition of this complex
\begin{eqnarray*}
\Omega^n(\geo PN\V/\U\geo)&=&\bigoplus_{p+q+r=n}\Omega^{p,q,r}(\geo PN\V/\U\geo)\\
&=&\bigoplus_{p+q_0+\cdots+q_p+r=n}\Omega^{p,q_0,\dots,q_p,r}(\geo PN\V/\U\geo),
\end{eqnarray*}
where $\Omega^{p,q_0,\dots,q_p,r}(\geo PN\V/\U\geo)$ is the set of forms of degree $p$ in the barycentric coordinates of the first simplex, of degree $q_0$ in the second and so on and finally of degree $r$ in some local coordinates on the nerve of the covering. This makes $\Omega^*(\geo PN\V/\U\geo)$ into a triple-complex. There is also a corresponding \v Cech-de Rham triple-complex
\[\check\Omega^{p,q,r}(\V/\U)=\bigoplus_{q_0+\dots+q_p=q}\Omega^r(P_pN\V/\U_{q_0,\dots, q_p})\]
with differentials
\begin{eqnarray*}
\partial':\check\Omega^{p,q,r}(\V/\U)\to \check\Omega^{p+1,q,r}(\V/\U)\\
\partial'':\check\Omega^{p,q,r}(\V/\U)\to \check\Omega^{p,q+1,r}(\V/\U)\\
\partial''':\check\Omega^{p,q,r}(\V/\U)\to \check\Omega^{p,q,r+1}(\V/\U)
\end{eqnarray*}
Here $\partial'=\sum(-1)^i \partial'_i$ where
\begin{equation*}
\partial'_i\alpha_{|j^0_0\dots j^{p+1}_{q_{p+1}}}=\left\{\begin{array}{ll}
0 &\textrm{if } q_i>0\\
\alpha_{|j^0_0\dots \hat{j^i_0}\dots j^{p+1}_{q_{p+1}}} & \textrm{if }q_i=0
\end{array}\right.
\end{equation*}
$\partial''$ and $\partial'''$ are usual \v Cech and de Rham differentials.

As in the simplicial case we have
\begin{proposition}
The map
\[I_\Delta:\Omega^{p,q,r}(\geo PN\V/\U\geo)\to \check\Omega^{p,q,r}(\V/\U)\]
given by 
\[I_\Delta(\omega)=\int_{\Delta^p\times \Delta^{q_0\dots q_p}}\omega_{q_0\dots q_p},\quad \textrm{for } \omega\in\Omega^{p,q_0,\dots,q_p,r}(\geo PN\V/\U\geo)\]
induces an isomorphism in cohomology. The right inverse is given on $\Delta^{k_0\dots k_p}\times P_pN\V/\U_{k_0\dots k_p}$ by
\[E(\omega)=p!q_0!\cdots q_p!\sum_{|J|=p}\sum_{|J_0|=q_0}\cdots\sum_{|J_p|=q_p} \omega_{J}\wedge\omega_{J_0}\wedge\cdots\wedge\omega_{J_p}\wedge d_{J_0\cdots J_p}^*\omega.\]
The $\omega_{J_j}$'s are the elementary forms on $\Delta^{q_j}$ and $d_{J_0\cdots J_p}$ are face maps as in the simplicial case.
\end{proposition}
\begin{proof}
The proof is the same as in the simplicial case (see e.g. \cite{D}).
\end{proof}
\begin{proposition}
The inclusion $\varepsilon:\sqcup V^i_j\to W_i$ induces the maps $\varepsilon_1^*$ and $\varepsilon_2^*$ in the following commutative diagram
\begin{eqnarray*}
\xymatrix{
\Omega^{p,r}(|N\W|) \ar[d]^{I_\Delta} \ar[r]^{\varepsilon_1^*} & \Omega^{p+r}(\geo PN\V/\U\geo) \ar[d]_{I_{\Delta}}\\
\check\Omega^{p,r}(\W)\ar[r]^{\varepsilon_2^*}& \check\Omega^{p+r}(\V/\U)
}
\end{eqnarray*}
They both induce isomorphisms in cohomology.
\end{proposition}
\begin{proof}
We first notice that both $\Omega^{p,0,r}(\geo PN\V/\U\geo)\cong\Omega^{p,r}(|N\V|)$ and $\check\Omega^{p,0,r}(\V/\U)\cong\Omega^{p,r}(\V)$ as double complexes, so in the following diagram
\begin{eqnarray*}
\xymatrix{
\Omega^{p,r}(|N\W|) \ar[d]^{I_\Delta} \ar[r]^{\varepsilon_1^*} & \Omega^{p,0,r}(\geo PN\V/\U\geo) \ar[d]_{I_{\Delta}}\\
\check\Omega^{p,r}(\W)\ar[r]^{\varepsilon_2^*}& \check\Omega^{p,0,r}(\V/\U)
}
\end{eqnarray*}
the $\varepsilon_i$'s are just refinement maps and thus induce isomorphisms in cohomology.

Now let us see that for fixed $p$ and $r$ the complex $\check\Omega^{p,q,r}(\V/\U)$ is exact, this will imply that
\[\check\Omega^{p,0,r}(\V/\U)\to\check\Omega^{p+r}(\V/\U)\]
is a cohomology isomorphism.

We first construct homomorphisms 
\[s^i:\check\Omega^{p,q_0,\dots, q_p,r}(\V/\U)\to\check\Omega^{p,q_0,\dots,q_i-1,\dots q_p}(\V/\U).\]
We choose partitions of unity on $W_i$ subordinate $\V^i=\{V^i_j\}_{j\in J_i}$ for each $i=i_0,\dots,i_p$ and set
\[s^i(\omega)_{j^0_0\dots j^i_{q_i-1}j^{i+1}_0\dots j^p_{q_p}}=(-1)^{q_0+\cdots+q_i}\sum_{j\in J_i}\phi^i_j\omega_{j^0_0\dots j^i_{q_i-1}jj^{i+1}_0\dots j^p_{q_p}}\]
Let $\tilde\delta^i=(-1)^{q^0+\cdots+q_{i-1}}\delta^i$, then we have
\[s^i\tilde\delta^j+\tilde\delta^js^i=0,\quad i\neq j\]
and
\[s^i\tilde\delta^i+\tilde\delta^is^i=\id\]
this gives
\[s^i\delta+\delta s^i=\id\]
for each $i$. So for fixed $p$ and $r$ the chain complex $\check\Omega^{q,r}(\V/\U)$ is exact.
\end{proof}
\begin{corollary}\label{corderham}
We have a quasi-isomorphism
\[\varepsilon'^*:\Omega^*(|N\W|)\to \Omega^*(\geo PN\V/\U\geo)\]
induced by the inclusion $\sqcup V^i_j\to W_i$.
\end{corollary}
\begin{rem}
The above result could have been obtained in a different manner. In the next section, we will construct a right inverse $\phi$ to $\varepsilon'$. We could then have constructed a homotopy $\phi\circ\varepsilon'\sim \id$ which would give us a chain homotopy directly on $\Omega^*(\geo PN\V/\U\geo)$.

\end{rem}


\section{Integration}\label{sectionintegration}
For a fiber bundle with compact oriented fibers, we want to define an integration map $\int:\Omega^{k+n}(|N\V|)\to \Omega^k(|N\U|)$ for coverings $\U$ and $\V$ coming from triangulations as in example \ref{eksprismaticnerve}. To do so, we define a map $|N\W|\to |N\V|$, and then our integration is given by pulling back forms by this map and then integrating along the fiber in $|N\W|\to |N\U|$. We define the map in two steps. First we have, similarly to the 'prismatic triangulation' map in example \ref{eksbundletriangulations}, a map $\lmap:\geo PN\V/\U\geo\to |N\V|$ defined on
\[\lmap:\Delta^p\times \Delta^{q_0\dots q_p}\times V_{j^0_0\dots j^p_{q_p}}\to\Delta^{p+q_0+\cdots+ q_p}\times V_{j^0_0\dots j^p_{q_p}}\]
by
\[\lmap(t,s^0,\dots,s^p,x)=(t_0s^0,\dots, t_ps^p,x)\]

Now recall that each $W_i$ is covered by $\V^i=\{\st(b^i_j)\}_{j\in J_i}$. Choose partitions of unity $\{\phi^i_j\}$ for $W_i$ subordinate $\V^i$ for each $i$. We are now ready to define 
\[\tilde\phi:|N\W|\to \geo PN\V/\U\geo\]
on $\Delta^p\times W_{i_0\dots i_p}$. Take $x\in W_{i_0\dots i_p}$ for each $i=i_0,\dots,i_p$ there is a minimal set $\{j^i_0,\dots j^i_{q_i}\}\in J_i$ so that 
\[\sum_{r=0}^{q_i} \phi^i_{j^i_r}(x)=1.\]
We then map
\[(t,x)\in \Delta^p\times W_{i_0\dots i_p}\]
to 
\[(t,\phi^{i_0}_{j^0_0}(x),\dots \phi^{i_0}_{j^0_{q_0}}(x),\dots,\phi^{i_p}_{j^p_{q_p}}(x),x)\in \Delta^p\times \Delta^{q_0\dots q_p}\times V_{j^0_0\cdots j^p_{q_p}}\]
\begin{rem}\label{remdimension}
Note that since the covering comes from a triangulation it has covering dimension $n+m$ so we have ensured that $q=\sum q_i\leq n$ for non-degenerate simplices.
\end{rem}
Now for $\omega\in \Omega^{n+k}(|N\V|)$ define
\[\left(\int_{[Y/Z]}\omega\right)_{|\Delta^p\times U_{i_0\dots i_p}}=\int_{\Delta^p\times W_{i_0\dots i_p}/\Delta^p\times U_{i_0\dots i_p}}\back\back\back\tilde\phi^* \lmap^*\omega,\]
where the right hand side denotes usual integration along the fibers.
\begin{theorem}
Given triangulations and partitions of unity as above, the following holds.
\paritem{$1.$}
Let $\omega\in\Omega^{*+n}(|N\V|)$ be a normal simplicial form, then $\int_{[Y/Z]}\omega$ is a well-defined normal simplicial form.
\paritem{$2.$} 
For $\omega\in\Omega^{*+n-1}(|N\V|)$ we have
\[\int_{[Y/Z]}d\omega=\int_{[\partial Y/Z]}\omega +(-1)^{n0}d\int_{[Y/Z]}\omega.\]
\end{theorem}
\begin{proof}
\paritem{$1.$} 
It is clear that $\int_{[Y/Z]}\omega$ is a well-defined simplicial form i.e. is compatible with respect to the degeneracy operators. Let us see that it is normal, that is
\[ (\eta_j\times \id)^*\left(\int_{[Y/Z]}\omega\right)^{(p)}=(\id\times s_j)^*\left(\int_{[Y/Z]}\omega\right)^{(p+1)}.\] 
We first notice that
\begin{eqnarray*}
(\eta_j\times \id)^*\left(\int_{[Y/Z]}\omega\right)_{|\Delta^p\times U_{i_0\dots i_p}} &=& (\eta_j\times \id)^*\int_{\Delta^p\times W_{i_0\dots i_p}/\Delta^p\times U_{i_0\dots i_p}}\back\back\back\tilde\phi^*\lmap^*\omega\\
&=& \int_{\Delta^{p+1}\times W_{i_0\dots i_p}/\Delta^{p+1}\times U_{i_0\dots i_p}}\back\back\back(\eta_j\times \id)^*\tilde\phi^*\lmap^*\omega\\
&=& \int_{\Delta^{p+1}\times W_{i_0\dots i_p}/\Delta^{p+1}\times U_{i_0\dots i_p}}\back\back\back\tilde\phi^*(\eta_j\times \id)^*\lmap^*\omega\\
&=& \int_{\Delta^{p+1}\times W_{i_0\dots i_p}/\Delta^{p+1}\times U_{i_0\dots i_p}}\back\back\back\tilde\phi^*(\lmap\circ(\eta_j\times \id))^*\omega
\end{eqnarray*}
and at the same time
\begin{eqnarray*}
(\id\times s_j)^*&&\back\back\left(\int_{[Y/Z]}\omega\right)_{|\Delta^{p+1}\times
 U_{i_0\dots i_ji_j\dots i_p}}\\ 
&=& (\id\times s_j)^*\int_{\Delta^{p+1}\times W_{i_0\dots i_ji_j\dots
 i_p}/\Delta^{p+1}\times U_{i_0\dots i_ji_j\dots i_p}}\back\back\back\tilde\phi^*\lmap^*\omega\\
&=& \int_{\Delta^{p+1}\times W_{i_0\dots i_j\dots i_p}/\Delta^{p+1}\times U_{i_0\dots i_j\dots i_p}}\back\back\back(\id \times s_j)^*\tilde\phi^*\lmap^*\omega\\
&=& \int_{\Delta^{p+1}\times W_{i_0\dots i_p}/\Delta^{p+1}\times U_{i_0\dots i_p}}\back\back\back\tilde\phi^*(\id\times s_j)^*\lmap^*\omega\\
&=& \int_{\Delta^{p+1}\times W_{i_0\dots i_p}/\Delta^{p+1}\times U_{i_0\dots i_p}}\back\back\back\tilde\phi^*(\lmap\circ(\id \times s_j))^*\omega.
\end{eqnarray*}
Hence we only need to show that $(\lmap\circ(\id\times s_j))^*\omega=(\lmap\circ(\eta_j\times \id))^*\omega$. This can be seen from the following commutative diagram
\begin{eqnarray*}
\xymatrix{
\Delta^p\times\Delta^{q_0\dots q_p}\times V_{j^0_0\dots j^p_{q_p}}  \ar[r]^{\lmap} & \Delta^{p+q}\times  V_{j^0_0\dots j^p_{q_p}}\\
\Delta^{p+1}\times\Delta^{q_0\dots q_p}\times V_{j^0_0\dots j^p_{q_p}} \ar[u]^{\eta_j\times \id}\ar[r]^{\tilde\lmap}\ar[d]^{\id\times s_j} & \Delta^{p+q+q_j+1}\times  V_{j^0_0\dots j^p_{q_p}}\ar[u]^{\tilde\eta\times \id}\ar[d]^{\tilde\sigma\circ (\id \times \tilde s)}\\
\Delta^{p+1}\times\Delta^{q_0\dots q_j q_j\dots q_p}\times V_{j^0_0\dots j^j_0\dots j^j_{q_j}j^j_0\dots j^p_{q_p}}  \ar[r]^{\lmap} & \Delta^{p+q+q_j+1}\times  V_{j^0_0\dots j^j_0\dots j^j_{q_j}j^j_0\dots j^p_{q_p}}
}
\end{eqnarray*}
where $q=\sum q_i$, $\tilde\lmap$ is given by
\begin{eqnarray*}
  \tilde\lmap(t,s^0,&\!\!\!\dots&\!\!\!,s^p,x)\\
  &=&(t_0s^0,\dots,t_js^j_0,t_{j+1}s^j_0,t_js^j_1,\dots, t_{j+1}s^j_{q_j},\dots,t_{p+1}s^p,x),
\end{eqnarray*}
$\tilde\eta$ (and similarly for $\tilde s$) is given by
\[\tilde\eta=\eta_{q_0+\cdots+q_{j-1}+j}\circ\eta_{q_0+\cdots+q_{j-1}+j+2}\circ\cdots\circ\eta_{q_0+\cdots+q_{j-1}+j+2q_j}\]
and finally $\tilde\sigma$ is the map that permutes the vertices in the simplex as in remark \ref{rempermutation}, so that by assumption $\tilde\sigma^*\omega=\omega$.
\paritem{$2.$}
Follows from the analogous formula for usual fiber integration.
\end{proof}

There is a map $\varepsilon':|N\V|\to |N\W|$ induced by the inclusions $\sqcup V^i_j\to W_i$ given on $\Delta^{p+q_0+\cdots+q_p}\times V_{j^0_0\dots j^p_{q_p}}$ by
\[\varepsilon'(t^0_0,\dots,t^0_{q_0},\dots,t^p_{q_p},x)=(\sum_j t^0_j,\dots,\sum_j t^p_j,x)\in \Delta^p\times W_{i_0\dots i_p}.\]
Since $\varepsilon'$ is left inverse to $\lmap\circ \tilde\phi$ the following lemma follows easily from the construction of the integral
\begin{lemma}\label{compatiblewithusualintegration}
The following diagrams commute
\begin{eqnarray*}
\begin{array}{cc}
\xymatrix{
\Omega^{*+n}(Y) \ar[d]_{\int_{Y/Z}} \ar[r]^{\varepsilon^*} & \Omega^{*+n}(|N\V|) \ar[d]^{\int_{[Y/Z]}}\\
\Omega^*(Z) \ar[r]^{\varepsilon^*} & \Omega^*(|N\U|)
}
&
\xymatrix{
\Omega^{*+n}(|N\W|) \ar[d]_{\int_{Y/Z}} \ar[r]^{\varepsilon'^*} & \Omega^{*+n}(|N\V|) \ar[dl]^{\int_{[Y/Z]}}\\
\Omega^*(|N\U|) &
}
\end{array}
\end{eqnarray*}
that is the integration of simplicial forms is compatible with the usual fiber integration.
\end{lemma}


\section{A combinatorial formula}
As in the case of a product bundle (cf. Gomi-Terashima \cite{GT}) there is also a combinatorial formula for the integration map in the case of a general fiber bundle. However, the resulting forms are only 'piece-wise' smooth. We start by introducing a new complex consisting of these forms.


\subsection{The triangulated nerve}
Given a triangulation $L$ of a smooth manifold $Z$ we have as mentioned earlier an open cover $\U$ given by the stars $\st(a)$ where $a$ is a point in $L^0$.

For every simplex $\sigma\in L$, the closed star 
\[\overline{\st(\sigma)}=\bigcup_{\tau\in L^n, \sigma\subseteq \tau} |\tau|\]
inherits a natural triangulation $L_{\sigma}$ from $L$. This gives a realisation $|L_\sigma|\cong \overline{\st(\sigma)}$. 
\begin{defi} The triangulated nerve $NL$ is the simplicial complex with $p$-simplices given by 
\[N_pL=\bigsqcup_{\sigma \in L^p} |L_\sigma|\]
and for $\sigma=[a_0,\dots,a_p]$ the face and degeneracy operators 
\[d_j:|L_{a_0\dots a_p}|\to |L_{a_0\dots \hat{a_j}\dots a_p}|\textrm{ and }s_j:|L_{a_0\dots a_p}|\to |L_{a_0\dots a_j a_j\dots a_p}|\]
are given by inclusions.
\end{defi}
 Our construction will give simplicial forms on $|NL|$. 

Recall that a form $\omega$ on a simplicial complex is a collection of forms $\omega=\{\omega^{(p)}\}
$ with $\omega^{(p)}\in\Omega^*(\Delta^p\times N_pL)$ satisfying the relation $(\varepsilon^j\times \id)^*\omega^{(p)}=(\id\times d_j)^*\omega^{(p-1)}$. But the $L_{\sigma}$'s, $\sigma\in L^p$, are simplicial sets too, so our forms $\omega^{(p)}$ actually live on 
\[\sqcup_{\sigma\in L^p} \sqcup_i\Delta^p\times \Delta^i\times L_{\sigma}^{(i)},\]
where $L^{(i)}_{\sigma}$ is the discrete set of $i$-simplices in $L_{\sigma}$.

Now much of what has been done in the previous sections carry over. We can define integral forms $\Omega^*_\Z(|NL|)\subseteq \Omega^*(|NL|)$ exactly as before and given triangulations of a fiber bundle as in example \ref{eksbundletriangulations} we also get triangulated nerves both of the base and the total space. We can also associate a prism complex to this situation in exactly the same way as in example \ref{eksprismaticnerve}. There is obviously also a map $\lmap:\geo PNK/L\geo\to |NK|$ as before.

Now let us show that with regard to cohomology it does not matter whether we use ordinary simplicial forms or simplicial forms on the triangulated nerves.

We introduce the simplicial manifold (with 'corners') $N\overline{\U}$ with
\[N_p\overline{\U}=\bigsqcup_{i_0,\dots,i_p}\overline{U_{i_0\dots i_p}}.\]
Since the cohomology of $\Omega^*(|N\U|)$ does not depend on the open cover and since forms on a closed subset are restrictions of forms on a larger open subset, we get that the restriction $\Omega^*(|N\overline{\U}|)\to \Omega^*(|N\U|)$ induces an isomorphism in cohomology.

\begin{proposition}
The map
\[\iota:\Omega^*(|N\bar\U|)\to \Omega^*(|NL|)\]
 induced by the homeomorphisms
\[|L_{a_0\dots a_p}|\cong \overline{\st([a_0,\dots,a_p])}\]
is an isomorphism in cohomology.
\end{proposition}
\begin{proof}
It follows readily from the following commutative diagram
\begin{eqnarray*}
\xymatrix{
\Omega^{p,q}(|N\bar \U|) \ar[r]^{\iota} \ar[d]^{I_\Delta} & \Omega^{p,q}(|NL|) \ar[d]_{I_{\Delta}}\\
\Omega^q(N_p\overline{\U}) \ar[r]^{\iota'}& \Omega^p(N_pL)
}
\end{eqnarray*}
since both vertical maps are isomorphisms in cohomology by the simplicial
de Rham theorem. In fact, the de Rham theorem also implies that the lower horisontal map induces an isomorphism in cohomology, since the map $\Omega^q(\overline{\st(\sigma)})\to \Omega^q(|L_{\sigma}|)$ is a cohomology isomorphism for all $\sigma\in L$.
\end{proof}

Now let us show that we can also represent a class in Deligne cohomology by a simplicial form on a triangulated nerve.

First, let $\U=\{U_i\}_{i\in I}$ be a covering of $Z$ and let $L$ be a triangulation, so that every closed star of $L$ lies inside an open set of $\U$. That is we have a map $\alpha:L^0\to I$ so that $\overline{\st(a)}\subseteq U_{\alpha(a)}$. This gives a chain map
\[T:\Omega^*(|N\U|)\to \Omega^*(|NL|)\]
\begin{proposition}
The map $T$ induces an isomorphism both in ordinary cohomology and between the cohomology of (\ref{simplicialcohom}) and
\begin{equation}\label{eqn:triangdeligne}
\Omega^{l-1}_{\R/\Z}(|NL|)\stackrel{d}{\to}\Omega^l_{\R/\Z}(|NL|)\stackrel{d}{\to}\Omega^{l+1}(|NL|)/\varepsilon^*\Omega^{l+1}(Z)
\end{equation}
\end{proposition}
\begin{proof}
Follows from the last proposition since $T$ is the composition of a refinement map and $\iota$.
\end{proof}
Hence smooth Deligne cohomology is represented by the cohomology of the sequence (\ref{eqn:triangdeligne}).


\subsection{The integration map}
We want to define an integration map 
\[\int_{K/L}:\Omega^{*+n}(|NK|)\to \Omega^*(|NL|).\]
First, for a simplex $\sigma=[a_0,\dots,a_p]\in L$ we define a map 
\begin{eqnarray*}
\aw:PC_{k,m}(K/L_\sigma)\to \bigoplus_{k_1+k_2=k} PC_{k_1,p}(K/\sigma)\otimes PC_{k_2,m}(K/L_\sigma).
\end{eqnarray*}
Let $(\tau,\eta)\in PC_{k,m}(K/L_\sigma)$ then since $\eta$ is a top-dimensional simplex in $L_\sigma$ we have $\sigma\subseteq \eta$ let $i_0,\dots, i_p\in \{0,\dots, n\}$ denote the indices of the corresponding vertices of $\sigma$ in $\eta$. Let us write $\tau$ as $\tau=[b_0^0,\dots b^0_{q_0}|\dots|b^m_0,\dots,b^m_{q_m}]$, where the $i$'th block, $|b^i_0,\dots,b^i_{q_i}|$, lies over the $i$'th vertex in $\eta$. For $0\leq s_j \leq q_{i_j}$ we define
\[\tau^{s_0\cdots s_p}=[b^{i_0}_0,\dots,b^{i_0}_{s_0}|\dots|b^{i_p}_0,\dots,b^{i_p}_{s_p}]\]
and
\[\tau_{s_0\cdots s_p}=[b^0_0,\dots,b^0_{q_0}|\dots|b^{i_j}_{s_j},\dots,b^{i_j}_{q_{i_j}}|\dots|b^m_0,\dots,b^m_{q_m}]\]
then our map is given by
\[\aw(\tau)=\sum_{0\leq s_j\leq q_{i_j}}\tau^{s_0\cdots s_p}\otimes\tau_{s_0\cdots s_p} .\]
That is an Alexander-Whitney type map with respect to each block of vertices in $\tau$ lying over a vertex in $\sigma$.

The following lemma is a straightforward computation similar to the proof of the usual $\aw$-map being a chain map. 
\begin{lemma}\label{ischainmap}
The map
\[\aw:PC_{k,m}(K/L_\sigma)\to \bigoplus_{k_1+k_2=k} PC_{k_1,p}(K/\sigma)\otimes PC_{k_2,m}(K/L_\sigma)\]
is a chain map with respect to the boundary map $\partial_F$ from example \ref{eksbundletriangulations}, that is
\[\aw \partial_F=\partial_F \aw.\]
\vspace{- 21 pt}
\begin{flushright} $ \square $ \end{flushright}
\end{lemma}
We have to specify $\int_{K/L}\omega\in \Omega^*(|NL|)$ as a form on $\Delta^p\times \eta$ for $\eta\in L_\sigma$. If for the moment we let $\eta$ be an $m$-simplex the formula is quite simple.

First pick an orientation of $\eta$, since the fibers of $\pi$ is oriented, this gives us an orientation of $Y_{|\eta}$ and hence a fundamental class $[Y_{|\eta}]\in PC_{n+m}(K/\sigma)$.

Now consider $NK_{|\pi^{-1}(|L_\sigma|)}$ as a subset of $|K|_{|\sigma}\times |K|_{||L_\sigma|}$. We will define $\int_{K/L}\omega_{|\Delta^p\times \eta}$ by restricting $\omega$ to $\aw([Y_{|\eta}])$ and integrate along the fiber over $\Delta^p\times \eta$.

Set $s=\sum_{i=0}^p s_i$, then our formula will be given by
\begin{equation}\label{combintformula}
\int_{K/L}\omega_{|\Delta^p\times\eta}=
\sum_{\tau\in PS_{n,m}(K/\eta)} \sum_{0\leq s_j\leq q_{i_j}}\varepsilon(\tau) \int_{\Delta^{p+s}\times \tau_{s_0\dots s_p}/\Delta^p\times\eta} \omega^{(p+s)}_{\tau^{s_0\dots s_p}},
\end{equation}
where $\omega^{(p+s)}_{\tau^{s_0\dots s_p}}\in \Omega^{*+n}(\Delta^p\times K_{\tau^{s_0\dots s_p}})$ and $\varepsilon(\tau)$ is the sign of $\tau$ in $[Y_{|\eta}]$. The integration shall be understood as follows: We restrict $\omega$ to $\Delta^{p+s}\times \tau_{s_0\dots s_p}$ and then integrate it along the fibers over $\Delta^p\times \eta$ with respect to the map $\Delta^{p+s}\to \Delta^p$ given by
\[(t_0,\dots t_{p+s})\mapsto(\sum_{i=0}^{s_0} t_i,\sum_{i=s_0+1}^{s_0+s_1+1} t_i,\dots,\sum_{i=s_0+\dots s_{p-1}+p}^{s_0+\dots s_p+p} t_i)\]
and the map $\tau_{s_0\dots s_p} \to \eta$ which is just the restriction of $\pi$.

\begin{rem}\label{prismviewokay}
In the above, we could also have chosen to pull $\omega$ back to $\Omega^{*+n}(\geo PNK/L_\sigma\geo)$ with $\lmap$ and then integrate with respect to the map
\[\Delta^p\times \Delta^{s_0 \dots s_p}\times \tau_{s_0\dots s_p}\to \Delta^p\times \eta\]
This gives the same result, but will be more convenient when we shall see that the two approaches to integration give the 'same' result.
\end{rem}

We still need to define the integral on $\Delta^p\times \eta'$ for $\eta'\in S_k(L_\sigma)$ a lower-dimensional simplex. This will actually just be the restriction of the integral on $\Delta^p\times \eta$ for $\eta$ a top-simplex such that $\eta'\subseteq \eta$, we shall see that this is independent of which top-simplex we choose (this also shows that the resulting form is really simplicial on $|L_\sigma|$).

Let us first take a look at what happens to the formula (\ref{combintformula}) when the integral is restricted to $\Delta^p\times \eta' \subseteq \Delta^p\times\eta$.

For a $\tau\in PS_{n,m}(K/\eta)$ we see that

\begin{equation}\label{partofintformula}
\int_{\Delta^{p+s}\times \tau_{s_0\dots s_p}/\Delta^p\times\eta} \omega^{(p+s)}_{\tau^{s_0\dots s_p}}
\end{equation}
restricted to $\Delta^p\times\eta'$ is non-zero exactly when $\tau_{s_0\dots s_p}\cap\pi^{-1}(\eta')$ and $\tau_{s_0\dots s_p}$ have the same dimension $r=n-s$ in the direction of the fiber. That is  $\tau_{s_0\dots s_p}\cap\pi^{-1}(\eta')\in PS_{r,k}(K/\eta')$ and $\tau_{s_0\dots s_p}\in PS_{r,m}(K/\eta)$ (the dimension in the direction of the fiber for a simplex in $K$ is given as $\dim_F \tau=\dim \tau-\dim \pi(\tau)$). 

Now in this case let $\alpha$ be the simplex in $L_\sigma$ 'spanned' by $\eta'$ and $\sigma$, then $\tau\cap\pi^{-1}(\alpha)$ is $n$-dimensional in the fiber direction, so over each $\tilde\eta \in S_m(L_\sigma)$, with $\eta'\subseteq\tilde\eta$, there is exactly one $\tilde\tau\in PS_{n,m}(K/\tilde\eta)$ with $\tau\cap\pi^{-1}(\alpha)\subseteq \tilde\tau$ and in the expression (\ref{partofintformula}) it would make no difference if we used $\tilde\tau$ instead of $\tau$.

We can also give an explicit formula in this case, but first we need some notation. For a top-simplex $\mu\in PS_{n,m}(K/L_\sigma)$ set $\tilde\mu=\mu\cap\pi^{-1}(\sigma)$. For a simplex  $\rho\in PS_{r,k}(K/\eta')$ let 
\begin{eqnarray*}
F\rho=\{\mu\in PS_{n,m}(K/L_\sigma) \!\!&\mid& \!\!\rho=\mu\cap\pi^{-1}(\eta'),\\
&&\dim_F \tilde\mu+\dim_F \rho-\dim_F(\tilde\mu\cap\rho)=n\}.
\end{eqnarray*}
Now write
\[\rho=[c^0_0,\dots,c^0_{r_0}|\dots|c^k_0,\dots,c^k_{r_k}]\textrm{ and } \tilde\mu=[b^0_0,\dots,b^0_{q_0}|\dots|b^p_0,\dots,b^p_{q_p}]\]
with $\mu\in F\rho$ and let $i_0,\dots,i_l\in\{0,\dots,p\}$ and $j_0,\dots, j_l\in \{0,\dots,k\}$ denote the coinciding blocks in $\tilde\mu$ and $\rho$, that is 
\[\tilde\mu\cap\rho=[b^{i_0}_0,\dots,b^{i_0}_{q_{i_0}}|\dots|b^{i_l}_0,\dots,b^{i_l}_{q_{i_l}}]=[c^{j_0}_0,\dots,c^{j_0}_{q_{j_0}}|\dots|c^{j_l}_0,\dots,c^{j_l}_{q_{j_l}}].\]
As before we set
\[\rho_{s_0\dots s_l}=[c^0_0,\dots,c^0_{r_0}|\dots|c^{j_\nu}_{s_\nu},\dots,c^{j_\nu}_{r_{j_\nu}}|\dots|c^k_0,\dots,c^k_{r_k}]\]
and
\[\tilde\mu^{s_0\dots s_l}=[b^0_0,\dots,b^0_{q_0}|\dots|b^{i_\nu}_0,\dots,b^{i_\nu}_{s_\nu}|\dots|b^p_0,\dots,b^p_{q_p}]\]
and then finally the integration formula is given on $\Delta^p\times \eta'$ by
\begin{equation}\label{formulasmallsimplex}
\sum_{\rho\in PS_{*,k}(K/\eta')} \sum_{\{\tilde\mu\mid \mu\in F\rho\}}\sum_{0\leq s_\nu\leq q_{i_\nu}} \varepsilon(\mu)\int_{\Delta^{p+s}\times \rho_{s_0\dots s_l}/\Delta^p\times\eta} \omega^{(p+s)}_{\tilde\mu^{s_0\dots s_l}}
\end{equation}

\begin{theorem}
\paritem{$1.$}
Let $\omega\in\Omega^{*+n}(|NK|)$ be a piece-wise smooth normal simplicial form; then $\int_{K/L}\omega$ is a well-defined piece-wise smooth normal simplicial form.
\paritem{$2.$} 
Let $\omega\in\Omega^{k+n-1}(|NK|)$, then we have a Stokes' theorem
\[\int_{K/L}d\omega=\int_{\partial_F K/L}\omega + (-1)^{n}d\int_{K/L}\omega\]
\paritem{$3.$} 
If $\partial Y=\emptyset$ then the map $\int_{K/L}:\Omega^{*+n}(|NK|)\to\Omega^*(|NL|)$ takes integral forms to integral forms and it induces a map $\pi_!:H_\del^{*+n}(Y,\Z)\to H_\del^*(Z,\Z)$ in smooth Deligne cohomology.
\end{theorem}
\begin{proof}
\paritem{$1.$}
This follows at once from the construction.
\paritem{$2.$}
First we observe that for $\omega\in \Omega^{k+n-1}(|NK|)$ we have on $\Delta^p\times \eta$ ($\eta\in L_\sigma^{(n)}$) 
\begin{eqnarray*} 
\int_{K/L}d\omega 
 & = &  \sum_{\tau\in PS_{n,m}(K/\eta)} \sum_{0\leq s_j\leq q_{i_j}} \int_{\Delta^{p+s}\times \tau_{s_0\dots s_p}/\Delta^p\times\eta} (d\omega)^{\tau^{s_0\dots s_p}}\\
 & = & \sum_{\tau\in PS_{n,m}(K/\eta)} \sum_{0\leq s_j\leq q_{i_j}} \int_{\partial_F(\Delta^{p+s}\times \tau_{s_0\dots s_p})/\Delta^p\times\eta} \omega^{(p+s)}_{\tau^{s_0\dots s_p}}+{}\\
 & & {} +  (-1)^{n} \sum_{\tau\in PS_{n,m}(K/\eta)} \sum_{0\leq s_j\leq q_{i_j}} d\int_{\Delta^{p+s}\times \tau_{s_0\dots s_p}/\Delta^p\times\eta} \omega^{(p+s)}_{\tau^{s_0\dots s_p}}\\
 & = & \sum_{\tau\in PS_{n,m}(K/\eta)} \sum_{0\leq s_j\leq q_{i_j}} \int_{\partial_F(\Delta^{p+s}\times \tau_{s_0\dots s_p})/\Delta^p\times\eta} \omega^{(p+s)}_{\tau^{s_0\dots s_p}}+{}\\
 & & {}+  (-1)^{n} d\int_{K/L}\omega.
\end{eqnarray*}
In this formula, we recognize the first terms as $\int_{\partial_F K/L}\omega$ since lemma \ref{ischainmap} gives us that $\partial_F AW([Y_{|\eta}])=AW(\partial_F[Y_{|\eta}])$.
Hence we have verified the formula for $\eta$ a top-dimensional simplex, and since the value of the integral on the other simplices is given by restrictions, the formula holds in general.
\paritem{$3.$}
If $\omega\in\Omega^{*+n}(|NK|)$ is integral, then we observe that the only non-zero terms in (\ref{combintformula}) are those for $s=n$, that is, the integration is only with respect to the map $\Delta^{p+n}\to\Delta^p$, and the resulting forms are then clearly integral. We also see that there is a result similar to lemma \ref{compatiblewithusualintegration}, so it is now clear that we have an induced map in Deligne cohomology.
\end{proof}

Now we are ready to compare the two integration maps. This comparison will also quite easily show that the first, smooth version of the integration map also takes integral forms to integral forms.

 First, choose a triangulation of the fiber bundle $Y\to Z$ and let $\V=\{V_j\}_{j\in J}$ and $\U=\{U_i\}_{i\in I}$ be the associated coverings by the stars. Now let $K$ and $L$ be subdivisions of these triangulations so that every closed star of $K$ and $L$ lies inside an open set of $\V$ and $\U$ respectively. Then we get maps 
\[T:\Omega^*(\geo PN\V/\U\geo)\to \Omega^*(\geo PNK/L\geo),\quad T':\Omega^*(|N\U|)\to \Omega^*(|NL|)\]
as above, inducing isomorphisms in cohomology.

\begin{lemma}\label{inttoint}
If $\partial Y=\emptyset$ then the map $\int_{[Y/Z]}:\Omega^{*+n}(|N\V|)\to\Omega^*(|N\U|)$ takes integral forms to integral forms and hence induces a map in smooth Deligne cohomology.
\end{lemma}
\begin{proof}
In the following, we will make use of remark \ref{prismviewokay}, that is, we will look at the integration map in terms of the prism complex.

Let $\beta\in\Omega^{k+n}(|N\V|)$ be an integral form. Now remark \ref{remdimension} ensures that the pull back $\lmap^*\beta\in\Omega^{k+n}(\geo PN\V/\U\geo)$ lies in the subcomplex $\bigoplus_{q\leq n}\Omega^{k+n-q,q,0}(\geo PN\V/\U\geo)$. Note also that everything besides the term in $\Omega^{k,n,0}(\geo PN\V/\U\geo)$ maps to zero under the integration map. Now the diagram
\begin{eqnarray}\label{compatiblediagram}
\xymatrix{
\Omega^{*+n}(|N\bar \W|) \ar@<-1ex>[r]_{\varepsilon'^*} \ar[dr]_{\int_{Y/Z}} & \Omega^{*+n}(\geo PN\V/\U\geo) \ar@<-1ex>[l]_{\tilde\phi^*}\ar[d]^{\int_{[Y/Z]}} \ar[r]^{T} & \Omega^{*+n}(\geo PNK/L\geo)\ar[d]^{\int_{K/L}} \\
 & \Omega^*(|N\U|) \ar[r]^{T'} & \Omega^*(|NL|)
}
\end{eqnarray}
where the commutativity of the triangle and the outer square  impliy that $\int_{K/L}T\varepsilon'^*=T'\int_{Y/Z}$. Furthermore, if we put $\beta'=\varepsilon'^*\tilde\phi^*\beta$ then by corollary \ref{corderham} we have 
\[\beta'-\beta=hd\beta+dh\beta,\] 
where $h$ is the homotopy operator inducing the chain homotopy
$\varepsilon'^*\tilde\phi^*\sim \id$. Note that since
$\tilde\phi\circ\varepsilon'$ is the identity in the variable of the first
simplex and in those on the nerve, $h$ maps $\Omega^{p,q,r}(\geo
PN\V/\U\geo )$ into 
\[\bigoplus_{q'<q}\Omega^{p,q-q'-1,r+q'}(\geo PN\V/\U\geo),\]
so the image of $h$ maps to zero under the integration map. Now by definition we have
\[\int_{[Y/Z]}\beta=\int_{Y/Z}\tilde\phi^*\beta=\int_{[Y/Z]}\beta'\]
and hence commutativity of the outer square in (\ref{compatiblediagram}) gives
\begin{eqnarray*}
T'\int_{Y/Z}\tilde\phi^*\beta&=&\int_{K/L}T\varepsilon'^*\tilde\phi^*\beta=\int_{K/L}T\beta'\\
&=&\int_{K/L}T\beta+\int_{K/L}T(hd\beta+dh\beta).
\end{eqnarray*}
Also, since $\int_{K/L}Tdh\beta=(-1)^{n-1}d\int_{K/L}Th\beta$ the last integral is zero. We therefore finally get
\[T'\int_{[Y/Z]}\beta=\int_{K/L}T\beta\]
and since the right side is clearly integral, as noted above, we conclude that $\int_{[Y/Z]}$ maps integral forms to integral forms. 
\end{proof}

\begin{theorem}\label{inducessamemap}
If $\partial Y=\emptyset$ then the maps
\[\int_{[Y/Z]}:\Omega^{*+n}(|N\V|)\to \Omega^*(|N\U|)\]
and 
\[\int_{K/L}:\Omega^{*+n}(|NK|)\to \Omega^*(|NL|)\]
induce the same map
\[\pi_!:H^{*+n}_\del(Y,\Z)\to H^*_\del(Z,\Z)\] 
 in smooth Deligne cohomology.
\end{theorem}
\begin{proof}
This is similar to the proof above. Taking $\omega\in\Omega^{*+n}(\geo PN\V/\U\geo)$ with $d\omega=\varepsilon^*\alpha-\beta$, we set $\omega'=\varepsilon'^*\tilde\phi^*\omega$ and get $\omega'-\omega=dh\omega+hd\omega$. As above we have
\begin{eqnarray*}
T'\int_{[Y/Z]}\omega&=&T'\int_{[Y/Z]}\omega'=\int_{K/L}T\omega'\\
&=&\int_{K/L}T\omega+d\int_{K/L}Th\omega+\int_{K/L}Thd\omega
\end{eqnarray*}
and, as in the proof of lemma \ref{inttoint} the last term vanishes, since $hd\omega=h\varepsilon^*\alpha-h\beta=-h\beta$ because $\varepsilon'^*\tilde\phi^*$ obviously acts as the identity on $\varepsilon^*\alpha$. Hence we get
\[T'\int_{[Y/Z]}\omega=\int_{K/L}T\omega+d\tau,\]
where $\tau=\int_{K/L}Th\omega$.
\end{proof}

\begin{corollary}
\[\pi_!:H^{*+n}_\del(Y,\Z)\to H^*_\del(Z,\Z)\]
is independent of choice of coverings, partition of unity and triangulations. In particular this proves theorem \ref{firstmain}.
\end{corollary}

Recall from example \ref{eksbundletriangulations} that in the case of a fiber bundle $Y \to Z$ with compact oriented fibers, where $\partial Y\neq \emptyset$ and a given triangulation of the bundle $\partial Y\to Z$, it is possible to extend this triangulation to a triangulation of the bundle $Y\to Z$. We shall see that the integral is independent of this extension, thus proving theorem \ref{secondmain}.
\begin{theorem}
Given a form $\omega\in \Omega^{*+n}(|N\V|)$ representing a class in
Deligne cohomology and a triangulation of $\partial Y\to Z$ compatible with the covering $\V$ and two extensions $|K_1|\to |L|$ and $|K_2|\to |L|$ of this to $Y\to Z$ then
\[\int_{K_1/L}T_1\omega \sim \int_{K_2/L}T_2\omega \quad \textrm{in }\Omega^*(|NL|),\]
where $T_i:\Omega^*(|N\V|)\to \Omega^*(|NK_i|)$, $i=1,2$ are given as above.
\end{theorem}
\begin{proof}
As in the proof of theorem \ref{inducessamemap} we have a $\omega'\in\Omega^{*+n}(\geo PN\V/\U\geo)$ so that $\omega'=\omega+dh\omega+hd\omega$, and we get, for $i=1,2$, 
\begin{eqnarray*}
T'\int_{[Y/Z]}\omega&=&T'\int_{[Y/Z]}\omega'=\int_{K_i/L}T_i\omega'\\
&=&\int_{K_i/L}T_i\omega+\int_{K_i/L}T_i(dh\omega+hd\omega)\\
&=&\int_{K_i/L}T_i\omega+\int_{K_i/L}T_idh\omega
\end{eqnarray*}
where $T':\Omega^*(|N\U|)\to \Omega^*(|NL|)$. Now the theorem follows from the fact that
\[\int_{K_i/L}T_idh\omega=d\int_{K_i/L}T_ih\omega \pm\int_{\partial_FK_i/L}T_ih\omega,\]
where the last term is easily seen to be independent of $i$.
\end{proof}


\section{Products}

The smooth Deligne cohomology groups comes with a product structure which has a quite simple description in the \v Cech-de Rham complex (see \cite{B} for a detailed description).

The usual wedge product on $\Omega^*(|N\U|)$ does not take integral forms to integral forms, so we have to define a new product with this property in order to get a product compatible with the usual product in the \v Cech-de Rham model.

First consider the maps
\[\pi_i:|P_1N\U|\to |N\U|,\quad i=1,2\]
where
\[\pi_1:\Delta^{q_0q_1}\times U_{i_0\dots i_{q_0+q_1+1}}\to \Delta^{q_0}\times U_{i_0\dots i_{q_0}}\]
is given by
\[(r^0,r^1,x)\mapsto(r^0,x)\]
and similarly 
\[\pi_2:\Delta^{q_0q_1}\times U_{i_0\dots i_{q_0+q_1+1}}\to \Delta^{q_1}\times U_{i_{q_0+1}\dots i_{q_0+q_1+1}}\]
is given by
\[(r^0,r^1,x)\mapsto (r^1,x).\]
For $t\in\Delta^1$ we have the map
\[\lmap_t:|P_1N\U|\to |N\U|\]
given by
\[\lmap_t(r^0,r^1,x)=(tr^0,(1-t)r^1,x).\]
It is clearly a homeomorphism for $t\in \stackrel{\circ}{\Delta^1}$.

The inverse is given as follows. Take an $(r_0,\dots,r_n,x)\in
\Delta^n\times U_{i_0\dots i_n}$ and choose $p$ so that $\sum_{i=0}^{p-1}r_i\leq t<\sum_{i=0}^p r_i$. Then
\begin{eqnarray*}
\lmap_t^{-1}(r_0,\dots,r_n,x)&=&\\
& &\back\back\back\back\back\back((\frac{r_0}{t},\dots,\frac{r_{p-1}}{t},1-\frac{\sum_{i=0}^{p-1}r_i}{t}),(1-\frac{\sum_{i={p+1}}^nr_i}{1-t},\frac{r_{p+1}}{1-t},\dots,\frac{r_n}{1-t}),s_px).
\end{eqnarray*}

Choose a smooth bump-function $\phi:\R\to\R$ so that the following holds
\paritem{$1.$} $\int_0^1\phi(t)dt=1$.
\paritem{$2.$} $\lim_{t\to 0}\phi(t)/t^p=0$ and $\lim_{t\to
  1}\phi(t)/(1-t)^p=0$ for all $p\in \N$.
We now have
\begin{defi}
The product
\[\wedge_1:\Omega^*(|N\U|)\times\Omega^*(|N\U|)\to \Omega^*(|N\U|).\]
is given by
\[\omega_1\wedge_1\omega_2:=\int_{\Delta^1}\phi(t) dt\wedge (\lmap_t^{-1})^*(\pi_1^*\omega_1\wedge\pi_2^*\omega_2)\]
\end{defi}
The choice of bump-function insures that there are no convergens problem. So this is well-defined and seen to give a normal simplicial form.

Most of the following proposition is trivial.
\begin{proposition}
\paritem{$1.$} Two different choices of bump-function gives chain-homotopic products.
\paritem{$2.$}
For $\omega_1\in\Omega^p(|N\U|)$ and $\omega_2\in\Omega^q(|N\U|)$ we have $d(\omega_1\wedge_1\omega_2)=d\omega_1\wedge_1\omega_2+(-1)^p\omega_1\wedge_1 d\omega_2$.
\paritem{$3.$} $I_\Delta:\Omega^*(|N\U|)\to \check\Omega^*(\U)$ is multiplicative.
\paritem{$4.$} If $\omega_1,\omega_2\in \Omega_\Z^*(|N\U|)$ then $\omega_1\wedge_1\omega_2\in\Omega_\Z^*(|N\U|)$.
\end{proposition}
\begin{proof}
\paritem{$1.$} This is trivial since two choices of bump-functions that satisfies the required conditions are certainly homotopic by a linear homotopy through such bump-functions.
\paritem{$2.$} Follows from the corresponding formula for the wedge product.
\paritem{$3.$} Suppose $\omega_1\wedge_1\omega_2\in\Omega^{n,m}(|N\U|)$ then we have
\begin{eqnarray*}
\int_{\Delta^n}(\omega_1\wedge_1\omega_2)_{i_0\dots i_n}&=&
\int_{\Delta^n}\int_{\Delta^1}\phi(t)dt\wedge(\lmap_t^{-1*}(\pi_1^*\omega_1\wedge\pi_2^*\omega_2))_{i_0\dots
  i_n}\\
&=&
\int_{\Delta^1\times\Delta^n}\back\phi(t)dt\wedge(\lmap_t^{-1*}(\pi_1^*\omega_1\wedge\pi_2^*\omega_2))_{i_0\dots
  i_n}\\
&=& \sum_{p+q=n}\int_{\Delta^1\times
  \Delta^p\times\Delta^q}\back\back\phi(t)dt\wedge\lmap_t^*\eta_p^*(\lmap_t^{-1*}(\pi_1^*\omega_1\wedge\pi_2^*\omega_2))_{i_0\dots i_n}\\
&=& \sum_{p+q=n}\int_{\Delta^1\times
  \Delta^p\times\Delta^q}\back\back\phi(t)dt\wedge(\omega_1)_{i_0\dots
  i_p}\wedge(\omega_2)_{i_p\dots i_n}\\
&=& \sum_{p+q=n}\int_{\Delta^p}(\omega_1)_{i_0\dots
  i_p}\wedge\int_{\Delta^q}(\omega_2)_{i_p\dots i_n}.
\end{eqnarray*}
So $I(\omega_1\wedge_1\omega_2)=I(\omega_1)\wedge I(\omega_2)$ as claimed.
\paritem{$4.$} This follows directly from the proof of $3.$
\end{proof}
\begin{rem}
Unfortunately the product is neither commutative nor associative on the chain level but $3$ above insures us that it is up to chain homotopy.
\end{rem}

The product is well-behaved with respect to the integration map in section \ref{sectionintegration}.
\begin{proposition}\label{intprod1}
For $\omega_1\in\Omega^{p+n}(|N\V|)$ and $\omega_2\in\Omega^{q}(|N\U|)$ we have
\begin{equation}\label{eqnintprod1}
\left(\int_{[Y/Z]}\omega_1\right)\wedge_1\omega_2=\int_{[Y/Z]}\omega_1\wedge_1\pi^*\omega_2.
\end{equation}
\end{proposition}
\begin{proof}
Note that
\[(\pi_1\times\pi_2)\circ \lmap_t^{-1}\circ\lmap\circ\tilde\phi=(\lmap\circ\tilde\phi\circ\pi_1\times\lmap\circ\tilde\phi\circ\pi_2)\circ \lmap_t^{-1}:|N\W|\to |N\V|\times |N\V|\]
So for a pair of forms $\omega,\tau\in\Omega^*(|N\U|)$ we get the relation
\[(\lmap\circ\tilde\phi)^*(\omega\wedge_1\tau)=(\lmap\circ\tilde\phi)^*\omega\wedge_1(\lmap\circ\tilde\phi)^*\tau.\]
This implies that
\begin{eqnarray*}
\int_{[Y/Z]}\omega_1\wedge_1\pi^*\omega_2 &=&
\int_{Y/Z}(\lmap\circ\tilde\phi)^*\omega_1\wedge_1(\lmap\circ\tilde\phi)^*\pi^*\omega_2\\
&=& \int_{Y/Z}(\lmap\circ\tilde\phi)^*\omega_1\wedge_1\pi^*\omega_2\\
&=& \left(\int_{Y/Z}(\lmap\circ\tilde\phi)^*\omega_1\right)\wedge_1\omega_2\\
&=& \left(\int_{[Y/Z]}\omega_1\right)\wedge_1\omega_2
\end{eqnarray*}
as stated above.
\end{proof}

Let us move on to Deligne cohomology where the product structure is a little different.
\begin{defi}
Let $\omega_1\in\Omega^p(|N\U|)$ and $\omega_2\in\Omega^q(|N\U|)$ be two forms representing classes in Deligne cohomology. That is $d\omega_i=\varepsilon^*\alpha_i-\beta_i$, where $\alpha_i$ is a global form and $\beta_i$ is integral. Then we define
\[\omega_1\tilde\wedge\omega_2:=\omega_1\wedge_1\varepsilon^*\alpha_2+(-1)^{p+1}\beta_1\wedge_1\omega_2\]
\end{defi}
Some calculations shows that this induces a well-defined product in Deligne cohomology. E.g. take another representative $\omega_1+\beta$ for the class $[\omega_1]$ then we have
\begin{eqnarray*}
(\omega_1+\beta)\tilde\wedge\omega_2 &=& (\omega_1+\beta)\wedge_1\varepsilon^*\alpha_2+(-1)^{p+1}(\beta_1+d\beta)\wedge_1\omega_2\\
&=& \omega_1\tilde\wedge\omega_2+\beta\wedge_1\varepsilon^*\alpha_2+(-1)^{p+1}d\beta\wedge_1\omega_2\\
&=& \omega_1\tilde\wedge\omega_2+d(\beta\wedge_1\omega_2)+\beta\wedge_1\beta_2\\
&\sim& \omega_1\tilde\wedge\omega_2.
\end{eqnarray*}
Notably $d(\omega_1\tilde\wedge\omega_2)=\varepsilon^*(\alpha_1\wedge_1\alpha_2)-\beta_1\wedge_1\beta_2$.

With this product on $\Omega^*(|N\U|)$, the map $I_\Delta$ of section \ref{afsnit:simplicielleformer} between the simplicial and the \v Cech-de Rham model for Deligne cohomology becomes an isomorphism of graded rings.

The next proposition follows directly from proposition \ref{intprod1}.
\begin{proposition}
For $\omega_1\in\Omega^{p+n}(|N\V|)$ and $\omega_2\in\Omega^{q}(|N\U|)$
representing classes in Deligne cohomology, we have
\begin{equation}
\left(\int_{[Y/Z]}\omega_1\right)\tilde\wedge\omega_2=\int_{[Y/Z]}\omega_1\tilde\wedge\pi^*\omega_2.
\end{equation}
\end{proposition}

\begin{rem}
The product described above simplifies proposition 5.17 in \cite {DK}. Since we can choose $\gamma_1\wedge_1\gamma_2$ as representative for $u_1\cup u_2$.
\end{rem}


\end{document}